\documentclass[12pt,english]{elsarticle}
\usepackage[T1]{fontenc}
\usepackage[latin9]{inputenc}
\usepackage[a4paper]{geometry}
\usepackage{color}
\usepackage{amsmath}
\usepackage{amsthm}
\usepackage{graphicx}
\usepackage{setspace}
\makeatletter
\theoremstyle{plain}
\newtheorem{thm}{\protect\theoremname}
\newenvironment{lyxlist}[1]
	{\begin{list}{}
		{\settowidth{\labelwidth}{#1}
		 \setlength{\leftmargin}{\labelwidth}
		 \addtolength{\leftmargin}{\labelsep}
		 }}
	{\end{list}}

\date{}

\makeatother

\usepackage{babel}
\providecommand{\theoremname}{Theorem}

\begin{document}

\title{\textbf{\huge{}Zhang-Zhang Polynomials of Ribbons}}

\author{\noindent \textbf{Bing-Hau He, Chien-Pin Chou, Johanna Langner, }\\
\textbf{and Henryk A. Witek{*}}}
\begin{onehalfspace}

\address{\noindent Department of Applied Chemistry and Institute of Molecular
Science,\\
 National Chiao Tung University, University Rd., 30010 Hsinchu, Taiwan}
\end{onehalfspace}

\address{\noindent {*} hwitek@mail.nctu.edu.tw}
\begin{abstract}
We report a closed-form formula for the Zhang-Zhang polynomial (\emph{aka}
ZZ polynomial or Clar covering polynomial) of an important class of
elementary pericondensed benzenoids $Rb\left(n_{1},n_{2},m_{1},m_{2}\right)$
usually referred to as ribbons. A straightforward derivation is based
on the recently developed interface theory of benzenoids {[}Langner
and Witek, \emph{MATCH Commun. Math. Comput. Chem.} \textbf{84}, 143\textendash 176
(2020){]}. The discovered formula provides compact expressions for
various topological invariants of $Rb\left(n_{1},n_{2},m_{1},m_{2}\right)$:
the number of Kekulé structures, the number of Clar covers, its Clar
number, and the number of Clar structures. The last two classes of
elementary benzenoids, for which closed-form ZZ polynomial formulas
remain to be found, are hexagonal flakes $O\left(k,m,n\right)$ and
oblate rectangles $Or\left(m,n\right)$. 
\end{abstract}
\maketitle

\section{Introduction\label{sec:Introduction}}

\begin{figure}
\centering{}\includegraphics[scale=0.6]{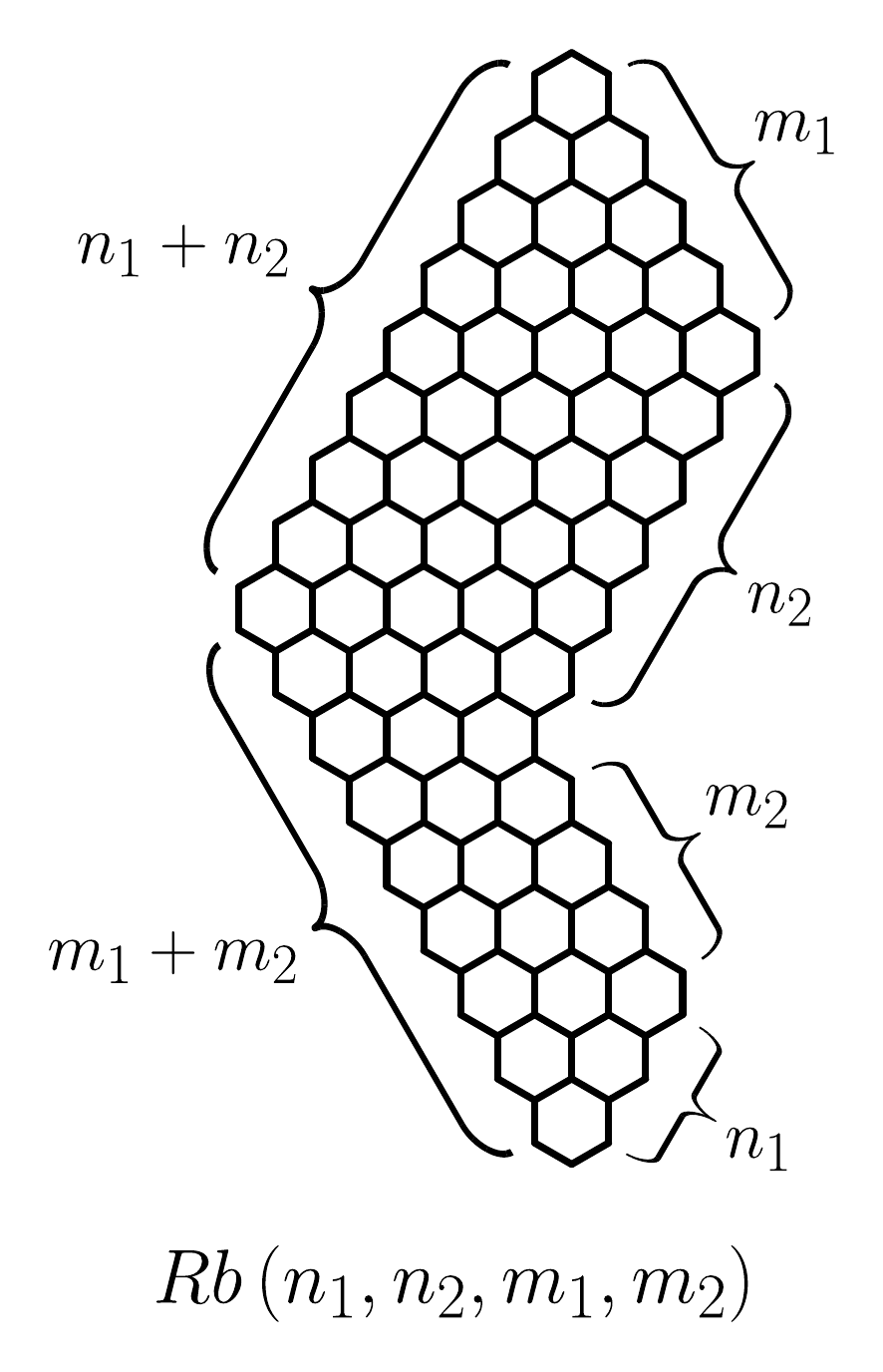}\caption{A graphical definition of ribbons $Rb\left(n_{1},n_{2},m_{1},m_{2}\right)$,
an important class of elementary pericondensed benzenoids. Here, $n_{1}=3$,
$n_{2}=6$, $m_{1}=5$ and $m_{2}=4$.\label{fig:graphdef}}
\end{figure}

Consider a regular pericondensed benzenoid depicted in Fig.~\ref{fig:graphdef},
 which can be fully characterized by specifying four structural parameters:
$n_{1}$, $n_{2}$, $m_{1}$, and $m_{2}$. Structures of this type
constitute an important family of elementary pericondensed benzenoids
and occupy a pronounced position in the general theory of Clar covers
among other, highly symmetric structures such as parallelograms $M\left(m,n\right)$
\citep{gordon1952theoryof,gutman2006zhangtextendash,chou2012analgorithm,chou2014closedtextendashform},
parallelogram chains \citep{he2020enum}, hexagons $O\left(k,m,n\right)$
\citep{gordon1952theoryof,ohkami1983topological,cyvin1986thenumber,chou2012zhangzhang,chou2014determination,he2020hexcorners},
oblate and prolate rectangles $Or\left(m,n\right)$ \citep{gutman2006algorithm,chou2012zhangzhang,chou2014determination}
and $Pr\left(m,n\right)$ \citep{zhang2000theclar,chou2012zhangzhang,chou2016closedform},
chevrons $Ch\left(k,m,n\right)$ \citep{gordon1952theoryof,cyvin1985enum,chou2012zhangzhang,chou2014closedtextendashform},
and generalized chevrons $Ch\left(k,m,n_{1},n_{2}\right)$ \citep{chou2014closedtextendashform}.
In situations when $m_{1}=n_{1}$, the structure shown in Fig.~\ref{fig:graphdef}
is traditionally referred to as a ribbon or a V-shaped benzenoid and
has been symbolically denoted by $V\left(k,m,n\right)$ \citep{cyvin1985enum,gutman1986topological,cyvin1988kekulestructures},
where $k=m_{1}=n_{1}$, $m=m_{1}+m_{2}$, and $n=n_{1}+n_{2}$. Here,
we consider a wider class of ribbons, allowing the parameters $m_{1}$
and $n_{1}$ to be different. We still refer to these structures as
ribbons, but we represent them by a new symbol $Rb\left(n_{1},n_{2},m_{1},m_{2}\right)$
capable of accommodating the extra new parameter in contrast to $V\left(k,m,n\right)$.
The Clar theory of generalized ribbons $Rb\left(n_{1},n_{2},m_{1},m_{2}\right)$\textemdash which
we attempt to construct in the current paper\textemdash clearly encompasses
the Clar theory of regular ribbons $V\left(k,m,n\right)$, similarly
as before we were able to show \citep{chou2014closedtextendashform}
that the Clar theory of generalized chevrons $Ch\left(k,m,n_{1},n_{2}\right)$
includes as a special case the Clar theory of regular chevrons $Ch\left(k,m,n\right)$.
The questions we want to answer in the current paper are: $(i)$~How
many Kekulé structures exist for $Rb\left(n_{1},n_{2},m_{1},m_{2}\right)$?
$(ii)$~How many Clar covers can be constructed for $Rb\left(n_{1},n_{2},m_{1},m_{2}\right)$?
$(iii)$~What is the Clar number of $Rb\left(n_{1},n_{2},m_{1},m_{2}\right)$?
$(iv)$~How many Clar structures can be constructed for $Rb\left(n_{1},n_{2},m_{1},m_{2}\right)$?
$(v)$~What is the ZZ polynomial of $Rb\left(n_{1},n_{2},m_{1},m_{2}\right)$?
Note that the solution to the last problem is sufficient for answering
all the posed here questions, justifying the title and the scope of
the current paper.

\section{\label{sec:Preliminaries}Preliminaries}

A benzenoid is a planar hydrocarbon $\boldsymbol{B}$ consisting of
fused benzene rings. From a graph theoretical point, $\boldsymbol{B}$
is defined as a 2-connected finite plane graph such that every interior
face is a regular hexagon \citep{gutman1989introduction}. A Kekulé
structure $\boldsymbol{K}$ is a resonance structure of $\boldsymbol{B}$
constructed using only double bonds \citep{kekule1866untersuchungen}.
A Clar cover $\boldsymbol{C}$ is a resonance structure of $\boldsymbol{B}$
constructed using double bonds and aromatic Clar sextets \citep{clar1972thearomatic}.
From a graph theoretical point, a Kekulé structure $\boldsymbol{K}$
is a spanning subgraph of $\boldsymbol{B}$ all of whose components
are $K_{2}$ and a Clar cover $\boldsymbol{C}$ is a spanning subgraph
of $\boldsymbol{B}$ whose components are either $K_{2}$ or hexagons
$C_{6}$. Note that most difficulties originating from this double
terminology can be circumvented if one establishes two correspondences:
$\text{a complete graph on 2 vertices }K_{2}\equiv\text{double bond}$
and a cycle graph on 6 vertices $C_{6}\equiv\text{aromatic Clar sextet}$.
The maximal number of hexagons $C_{6}$ that can be accommodated in
$\boldsymbol{C}$ is referred to as the Clar number $Cl$ of $\boldsymbol{B}$
\citep{clar1972thearomatic,gutman1985clarformulas}. The Clar covers
with $Cl$ aromatic sextets $C_{6}$ are referred to as the Clar structures
of $\boldsymbol{B}$ \citep{clar1972thearomatic,gutman1985clarformulas}.
The Clar covers with $k$ aromatic sextets $C_{6}$ are referred to
as the Clar covers of order $k$. If we represent the number of Clar
covers of order $k$ for $\boldsymbol{B}$ by $c_{k}$, we can define
a combinatorial polynomial
\begin{equation}
\text{ZZ}\left(\boldsymbol{B},x\right)=\sum_{k=0}^{Cl}c_{k}\,x^{k}\label{eq:zz-1}
\end{equation}
usually referred to as the Clar covering polynomial of $\boldsymbol{B}$
or the Zhang-Zhang polynomial of $\boldsymbol{B}$ or, shortly, the
ZZ polynomial of $\boldsymbol{B}$. \citep{zhang1996theclar,zhang1996theclar2,zhang1997theclar,zhang2000theclar}
Clearly, the ZZ polynomial of $\boldsymbol{B}$ has the following
inviting properties:
\begin{itemize}
\item The number of Kekulé structures of $\boldsymbol{B}$ is given by $K\left\{ \boldsymbol{B}\right\} =c_{0}=\text{ZZ}\left(\boldsymbol{B},0\right)$.
\item The number of Clar covers of $\boldsymbol{B}$ is given by $C\left\{ \boldsymbol{B}\right\} =c_{0}+\cdots+c_{Cl}=\text{ZZ}\left(\boldsymbol{B},1\right)$.
\item The Clar number of $\boldsymbol{B}$ is given by $Cl=\deg\left(\text{ZZ}\left(\boldsymbol{B},x\right)\right)$.
\item The number of Clar structures of $\boldsymbol{B}$ is given by $c_{Cl}=\text{coeff }\left(\text{ZZ}\left(\boldsymbol{B},x\right),x^{Cl}\right)=\left(Cl\right)!\frac{d}{dx^{Cl}}\text{ZZ}\left(\boldsymbol{B},x\right)$.
\end{itemize}
\noindent These relations demonstrate the claim made at the end of
Section~\ref{sec:Introduction}, where we have written that determination
of the ZZ polynomial of $\boldsymbol{B}$ answers most graph-theoretically
relevant questions about $\boldsymbol{B}$. Zhang and Zhang were able
to demonstrate that the ZZ polynomials possess a rich structure of
recursive decomposition properties \citep{zhang1996theclar,zhang1996theclar2,zhang2000theclar},
which enable their fast and robust computations in practical applications.
(See for example \textbf{Properties 1\textendash 7} in \citep{chou2012analgorithm}.)
Consequently, the ZZ polynomial of an arbitrary benzenoid $\boldsymbol{B}$
can be efficiently computed using recursive decomposition algorithms
\citep{gutman2006algorithm,chou2012analgorithm,chou2014zzdecomposer}
or determined using interface theory of benzenoids \citep{langner2018algorithm3,langner2017zigzagConnectivityGraph2,langner2019IFTTheorems5,langner2019BasicApplications6}.
A useful practical tool for determination of ZZ polynomials is ZZDecomposer
\citep{chou2014zzdecomposer,chou2014determination}. With this freely
downloadable \citep{ZZDecomposerDownload1,ZZDecomposerDownload2}
software, one can conveniently define a graph representation corresponding
to a given benzenoid $\boldsymbol{B}$ using a mouse drawing pad and
subsequently use it to find the ZZ polynomial of $\boldsymbol{B}$,
generate the set of Clar covers of $\boldsymbol{B}$, and determine
its structural similarity to other, related benzenoids. 

In the most typical depth-decomposition mode of ZZDecomposer, used
below in Fig.~\ref{fig:rb1dec} to prove Eqs.~(\ref{eq:ZZRb1})
and~(\ref{eq:ZZRb2}), ZZDecomposer generates a recurrence relation
for the analyzed benzenoid structure, which relates its ZZ polynomial
to the ZZ polynomials of structurally related benzenoids and often
allows for determination of a closed-form formulas for the whole family
of structurally similar benzenoids. Another useful feature of ZZDecomposer
is generating vector graphics that can be easily incorporated in publications. 

\section{Heuristic determination of the ZZ polynomial from recurrence relations}

\begin{figure}
\centering{}\includegraphics[clip,scale=0.14]{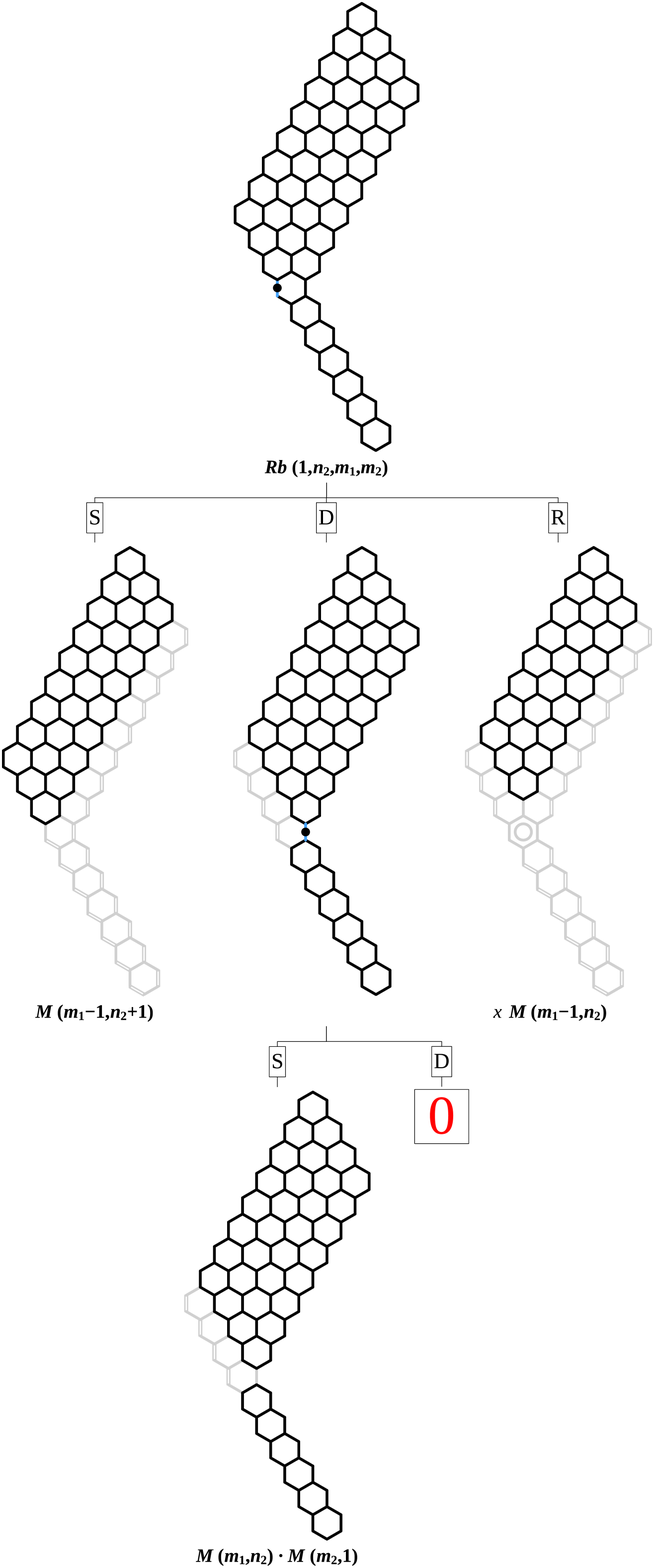}\includegraphics[clip,scale=0.14]{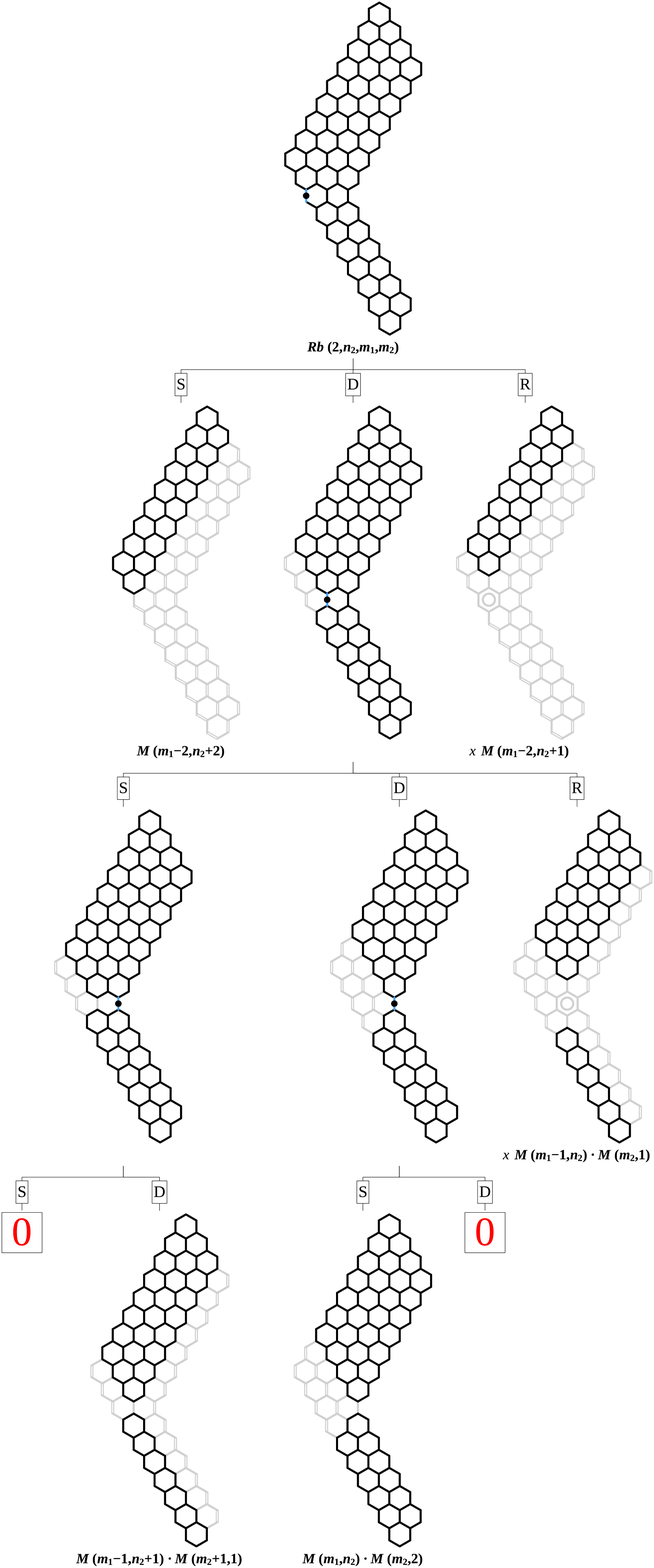}\caption{Multi-step recursive decomposition of the ribbons $Rb\left(1,n_{2},m_{1},m_{2}\right)$
and $Rb\left(2,n_{2},m_{1},m_{2}\right)$ with respect to the blue
edges (marked with a black circle in the middle) allows us to represent
their ZZ polynomials as a function of ZZ polynomials of various parallelograms
$M\left(m,n\right)$ in the form of the recurrence relation given
in Eqs.~(\ref{eq:ZZRb1}) and~(\ref{eq:ZZRb2}). The symbol $\boxed{{\color{red}\boldsymbol{0}}}$
denotes a non-Kekuléan decomposition pathway. Here, $m_{1}=4$, $m_{2}=6$
and $n_{2}=8$.\label{fig:rb1dec}}
\end{figure}

Before presenting a formal derivation of the ZZ polynomial for the
ribbon $Rb\left(n_{1},n_{2},m_{1},m_{2}\right)$ in Section~\ref{sec:IT},
first we discuss here a heuristic reasoning suggesting how such formulas
can be discovered for a general benzenoid $\boldsymbol{B}$. This
goal can be readily achieved using ZZDecomposer described in Section~\ref{sec:Preliminaries}
by considering the first two members of this family of structures,
$Rb\left(1,n_{2},m_{1},m_{2}\right)$ and $Rb\left(2,n_{2},m_{1},m_{2}\right)$,
and performing their multi-step recursive decompositions with respect
to the covering character of the vertical edges depicted in blue (with
a black dot) in Fig.~\ref{fig:rb1dec}.
 The process of assigning to these edges single bond covering $\boxed{\text{S}}$,
double bond covering $\boxed{\text{D}}$, or aromatic ring covering
$\boxed{\text{R}}$, as it is demonstrated in Fig.~\ref{fig:rb1dec},
allows us to partition the set of Clar covers of $Rb\left(1,n_{2},m_{1},m_{2}\right)$
and $Rb\left(2,n_{2},m_{1},m_{2}\right)$ into three and five, respectively,
subsets, each of them consisting of a region of fixed bonds separating
parallelogram-shaped regions with not fixed bonds. Consequently, the
ZZ polynomials of $Rb\left(1,n_{2},m_{1},m_{2}\right)$ and $Rb\left(2,n_{2},m_{1},m_{2}\right)$
can be expressed by summing over the subsets ZZ polynomials. Moreover,
the formulas can be written compactly in terms of the ZZ polynomials
of the parallelograms $M\left(m,n\right)$ as
\begin{eqnarray}
\text{ZZ}\left(Rb\left(1,n_{2},m_{1},m_{2}\right),x\right) & = & \text{ZZ}\left(M\left(m_{1},n_{2}\right),x\right)\cdot\text{ZZ}\left(M\left(m_{2},1\right),x\right)\label{eq:ZZRb1}\\
 & + & \text{ZZ}\left(M\left(m_{1}-1,n_{2}+1\right),x\right)+x~\text{ZZ}\left(M\left(m_{1}-1,n_{2}\right),x\right)\nonumber 
\end{eqnarray}
and
\begin{eqnarray}
\text{ZZ}\left(Rb\left(2,n_{2},m_{1},m_{2}\right),x\right) & = & \text{ZZ}\left(M\left(m_{1},n_{2}\right),x\right)\cdot\text{ZZ}\left(M\left(m_{2},2\right),x\right)\label{eq:ZZRb2}\\
 & + & \text{ZZ}\left(M\left(m_{1}-1,n_{2}+1\right),x\right)\cdot\text{ZZ}\left(M\left(m_{2}+1,1\right),x\right)\nonumber \\
 & + & x~\text{ZZ}\left(M\left(m_{1}-1,n_{2}\right),x\right)\cdot\text{ZZ}\left(M\left(m_{2},1\right),x\right)\nonumber \\
 & + & \text{ZZ}\left(M\left(m_{1}-2,n_{2}+2\right),x\right)+x~\text{ZZ}\left(M\left(m_{1}-2,n_{2}\right),x\right)\nonumber 
\end{eqnarray}
It is easy to notice regularities and patterns in these formulas.
Remembering that $\text{ZZ}\left(M\left(m,0\right),x\right)=1$, we
can rewrite Eqs.~(\ref{eq:ZZRb1}) and~(\ref{eq:ZZRb2}) for situations
when $n_{1}\leq m_{1}$ in the following form
\begin{equation}
\text{ZZ}\left(Rb\left(n_{1},n_{2},m_{1},m_{2}\right),x\right)=~~~~~~~~~~~~~~~~~~~~~~~~~~~~~~~~~~~~~~~~~~~~~~~~~~~~~~~~~~~~~~~~~~~~\label{eq:ZZRbM0}
\end{equation}
\[
~~~~~=~\sum_{k=0}^{n_{1}}\text{ZZ}\left(M\left(m_{1}-k,n_{2}+k\right),x\right)\cdot\text{ZZ}\left(M\left(m_{2}+k,n_{1}-k\right),x\right)+
\]
\[
~~~~~~~~~~~~~~~+~x~\sum_{k=1}^{n_{1}}\text{ZZ}\left(M\left(m_{1}-k,n_{2}-1+k\right),x\right)\cdot\text{ZZ}\left(M\left(m_{2}-1+k,n_{1}-k\right),x\right)
\]
The final generalization applies to situations when $n_{1}>m_{1}$.
It is clear that the ZZ polynomials of $Rb\left(n_{1},n_{2},m_{1},m_{2}\right)$
and $Rb\left(m_{1},m_{2},n_{1},n_{2}\right)$ should be identical
as there is a clear (horizontal mirror reflection) isomorphism between
the sets of Clar covers of both structures. Indeed, Eq.~(\ref{eq:ZZRbM0})
reflects this symmetry except for the upper summations limits; it
is easy to see that the appropriate change relies on replacing $n_{1}$
by $\min\left(n_{1},m_{1}\right)$. Consequently, the general formula
for the ZZ polynomial of $Rb\left(n_{1},n_{2},m_{1},m_{2}\right)$
must have the following symmetric form
\begin{equation}
\text{ZZ}\left(Rb\left(n_{1},n_{2},m_{1},m_{2}\right),x\right)=~~~~~~~~~~~~~~~~~~~~~~~~~~~~~~~~~~~~~~~~~~~~~~~~~~~~~~~~~~~~~~~~~~~~\label{eq:ZZRbM}
\end{equation}
\begin{equation}
~~~=~\sum_{k=0}^{\min\left(n_{1},m_{1}\right)}\text{ZZ}\left(M\left(m_{1}-k,n_{2}+k\right),x\right)\cdot\text{ZZ}\left(M\left(m_{2}+k,n_{1}-k\right),x\right)+\label{eq:sum1}
\end{equation}
\begin{equation}
~~~~~~~~~~~~~~~~+~x~\sum_{k=1}^{\min\left(n_{1},m_{1}\right)}\text{ZZ}\left(M\left(m_{1}-k,n_{2}-1+k\right),x\right)\cdot\text{ZZ}\left(M\left(m_{2}-1+k,n_{1}-k\right),x\right)\label{eq:sum2}
\end{equation}
This formula can be further simplified by substituting an explicit
form of the ZZ polynomial for the parallelogram $M\left(m,n\right)$
\begin{equation}
\text{ZZ}\left(M\left(m,n\right),x\right)={}_{2}F_{1}\left[\begin{array}{c}
{\scriptstyle -m,-n}\\
{\scriptstyle 1}
\end{array};1+x\right]=\sum_{j=0}^{\min\left(m,n\right)}\binom{m}{j}\binom{n}{j}\left(1+x\right)^{j}\label{eq:ZZMmn}
\end{equation}
Introducing this formula into Eq.~(\ref{eq:ZZRbM}), we obtain
\begin{equation}
\text{ZZ}\left(Rb\left(n_{1},n_{2},m_{1},m_{2}\right),x\right)=~~~~~~~~~~~~~~~~~~~~~~~~~~~~~~~~~~~~~~~~~~~~~~~~~~~~~~~~~~~~~~~~~~~~\label{eq:ZZRb}
\end{equation}
\[
~~~=~\sum_{k=0}^{\min\left(n_{1},m_{1}\right)}{}_{2}F_{1}\left[\begin{array}{c}
{\scriptstyle -m_{1}+k,-n_{2}-k}\\
{\scriptstyle 1}
\end{array};1+x\right]~{}_{2}F_{1}\left[\begin{array}{c}
{\scriptstyle -n_{1}+k,-m_{2}-k}\\
{\scriptstyle 1}
\end{array};1+x\right]+
\]
\[
~~~~~~~~+~x\sum_{k=1}^{\min\left(n_{1},m_{1}\right)}{}_{2}F_{1}\left[\begin{array}{c}
{\scriptstyle -m_{1}+k,-n_{2}+1-k}\\
{\scriptstyle 1}
\end{array};1+x\right]~{}_{2}F_{1}\left[\begin{array}{c}
{\scriptstyle -n_{1}+k,-m_{2}+1-k}\\
{\scriptstyle 1}
\end{array};1+x\right]
\]
or
\begin{equation}
\text{ZZ}\left(Rb\left(n_{1},n_{2},m_{1},m_{2}\right),x\right)=~~~~~~~~~~~~~~~~~~~~~~~~~~~~~~~~~~~~~~~~~~~~~~~~~~~~~~~~~~~~~~~~~~~~\label{eq:ZZRb-1}
\end{equation}
\[
=~\sum_{k=0}^{\min\left(n_{1},m_{1}\right)}\sum_{j=0}^{m_{1}-k}\sum_{i=0}^{n_{1}-k}\binom{m_{1}-k}{j}\binom{n_{2}+k}{j}\binom{m_{2}+k}{i}\binom{n_{1}-k}{i}\left(1+x\right)^{i+j}+
\]
\[
~~~~~~~~+~x\sum_{k=1}^{\min\left(n_{1},m_{1}\right)}\sum_{j=0}^{m_{1}-k}\sum_{i=0}^{n_{1}-k}\binom{m_{1}-k}{j}\binom{n_{2}-1+k}{j}\binom{m_{2}-1+k}{i}\binom{n_{1}-k}{i}\left(1+x\right)^{i+j}
\]
Numerical experiments performed with ZZDecomposer for various values
of the parameters $n_{1}$, $n_{2}$, $m_{1}$, and $m_{2}$ show
that formulas given by Eqs.~(\ref{eq:ZZRb}) and~(\ref{eq:ZZRb-1})
are indeed correct. Formal demonstration of correctness of Eqs.~(\ref{eq:ZZRb})
and~(\ref{eq:ZZRb-1}) is presented in the next Section; the proof
is based on the recently developed interface theory of benzenoids
\citep{langner2019IFTTheorems5,langner2019BasicApplications6}.

As we mentioned earlier, the ribbon $Rb\left(n_{1},n_{2},m_{1},m_{2}\right)$
has been also denoted in the earlier literature by the symbol $V\left(k,m,n\right)$,
where $k=m_{1}=n_{1}$, $m=m_{1}+m_{2}$, and $n=n_{1}+n_{2}$ or
by $V\left(k_{1},k_{2},m,n\right)$, where $k_{1}=n_{1}$, $k_{2}=m_{1}$,
$m=m_{1}+m_{2}$, and $n=n_{1}+n_{2}$. Therefore, for consistency,
we also give explicit formulas for the ZZ polynomials of $V\left(k,m,n\right)$
and $V\left(k_{1},k_{2},m,n\right)$ using their structural parameters
in the formulas. We have
\begin{equation}
\text{ZZ}\left(V\left(k,m,n\right),x\right)=~~~~~~~~~~~~~~~~~~~~~~~~~~~~~~~~~~~~~~~~~~~~~~~~~~~~~~~~~~~~~~~~~~~~\label{eq:ZZRb-1-1}
\end{equation}
\[
=~\sum_{s=0}^{k}\sum_{j=0}^{k-s}\sum_{i=0}^{k-s}\binom{k-s}{j}\binom{n-k+s}{j}\binom{m-k+s}{i}\binom{k-s}{i}\left(1+x\right)^{i+j}+
\]
\[
~~~~~~~~+~x\sum_{s=1}^{k}\sum_{j=0}^{k-s}\sum_{i=0}^{k-s}\binom{k-s}{j}\binom{n-k-1+s}{j}\binom{m-k-1+s}{i}\binom{k-s}{i}\left(1+x\right)^{i+j}
\]
and
\begin{equation}
\text{ZZ}\left(V\left(k_{1},k_{2},m,n\right),x\right)=~~~~~~~~~~~~~~~~~~~~~~~~~~~~~~~~~~~~~~~~~~~~~~~~~~~~~~~~~~~~~~~~~~~~\label{eq:ZZRb-1-2}
\end{equation}
\[
=~\sum_{k=0}^{\min\left(k_{1},k_{2}\right)}\sum_{j=0}^{k_{2}-k}\sum_{i=0}^{k_{1}-k}\binom{k_{2}-k}{j}\binom{n-k_{1}+k}{j}\binom{m-k_{2}+k}{i}\binom{k_{1}-k}{i}\left(1+x\right)^{i+j}+
\]
\[
~~~~~~~~+~x\sum_{k=1}^{\min\left(k_{1},k_{2}\right)}\sum_{j=0}^{k_{2}-k}\sum_{i=0}^{k_{1}-k}\binom{k_{2}-k}{j}\binom{n-k_{1}-1+k}{j}\binom{m-k_{2}-1+k}{i}\binom{k_{1}-k}{i}\left(1+x\right)^{i+j}.
\]
Note that analogous formulas for the number of Kekulé structures given
by Cyvin and Gutman (as Eq.~(19) of \citep{gutman1986topological})
have both of the inner summations evaluated to a closed binomial form.
This is indeed possible for Kekulé structures, for which $x=0$, where
the following binomial identity (Eq. (5.22) of \citep{gould2015combi})
can be used.
\begin{equation}
\sum_{j=0}^{b}\binom{b}{j}\binom{v-b}{j}=\binom{v}{b}\label{eq:binid}
\end{equation}
In the case of the ZZ polynomial, Eq.~(\ref{eq:binid}) takes on
the following form
\[
\sum_{j=0}^{b}\binom{b}{j}\binom{v-b}{j}\left(1+x\right)^{j}={}_{2}F_{1}\left[\begin{array}{c}
{\scriptstyle -b,-v+b}\\
{\scriptstyle 1}
\end{array};1+x\right]
\]
which renders Eqs.~(\ref{eq:ZZRb-1-1}) and (\ref{eq:ZZRb-1-2})
in the hypergeometric form analogous to Eq.~(\ref{eq:ZZRb}). These
hypergeometric functions reduce to the obvious polynomial form or
to Jacobi polynomials, and no other functional identities exist that
would allow to express them as some well-known functions.

\section{\label{sec:IT}Formal derivation of the ZZ polynomial from the interface
theory of benzenoids}

\begin{figure}
\centering{}\includegraphics[scale=0.6]{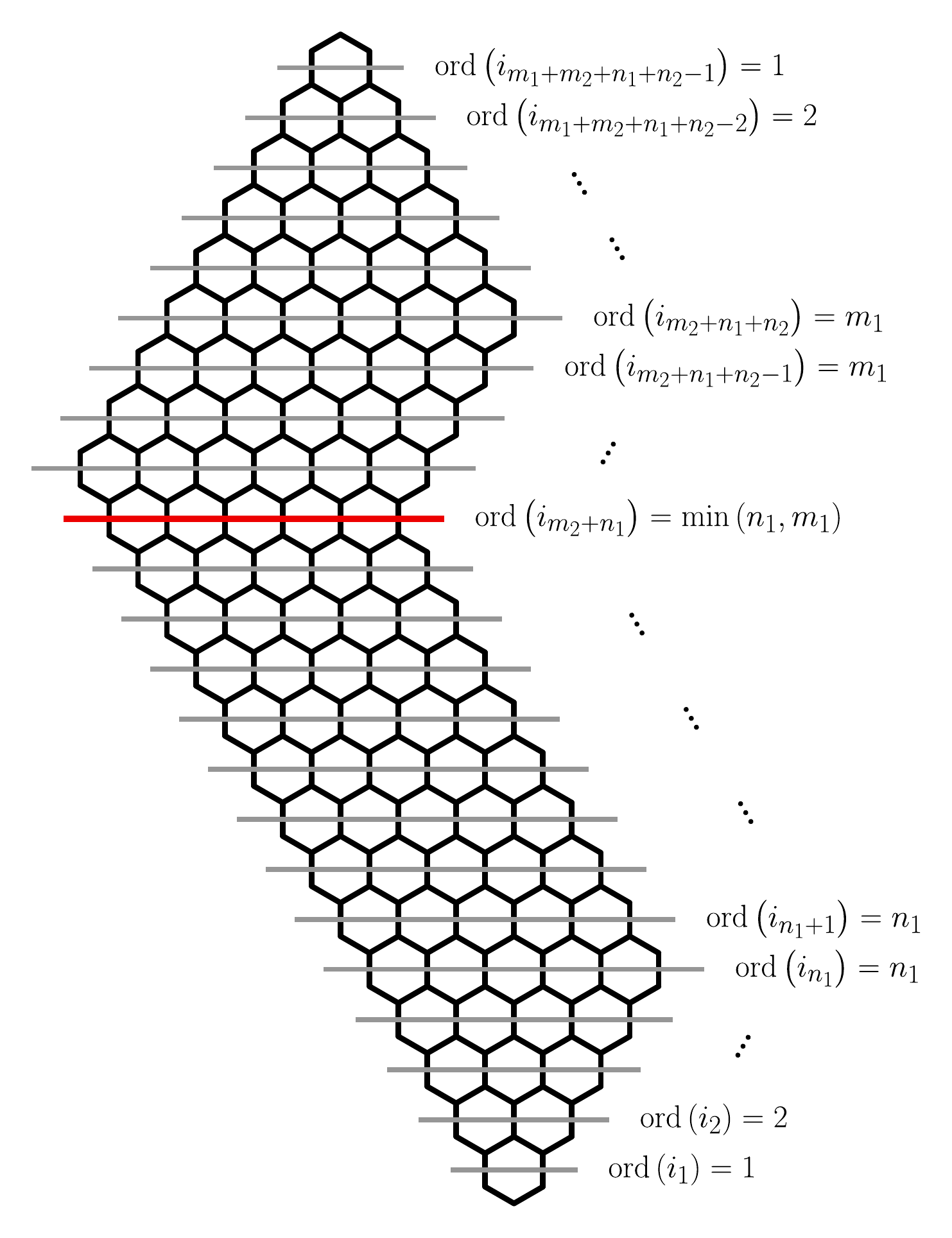}\caption{A system of elementary cuts (represented by horizontal lines) for
the ribbon $Rb\left(n_{1},n_{2},m_{1},m_{2}\right)$ defines $n_{1}+m_{1}+m_{2}+n_{2}-1$
interfaces $i_{1},\ldots,i_{n_{1}+m_{2}+n_{2}+m_{1}-1}$. Orders of
these interfaces can be calculated according to the interfaces theory.
The central interface $i_{n_{1}+m_{2}},$depicted using red horizontal
line, has an order equal to the minimum of the two structural constants
of the ribbon, $m_{1}$ and $n_{1}$. Here, $n_{1}=5$, $n_{2}=4$,
$m_{1}=6$ and $m_{2}=9$.\label{fig:interfaces}}
\end{figure}

Consider the ribbon $\boldsymbol{B}\equiv Rb\left(n_{1},n_{2},m_{1},m_{2}\right)$
in the orientation shown in Fig.~\ref{fig:graphdef}. We introduce
a system of $m_{1}+n_{2}+m_{2}+n_{1}-1$ elementary cuts $I_{k}$
intersecting the vertical edges of $\boldsymbol{B}$ in the way shown
in Fig.~\ref{fig:interfaces}. The set of vertical edges intersected
by the elementary cut $I_{k}$ is referred to as the \emph{interface}
$i_{k}$ of $\boldsymbol{B}$. It is convenient to augment the set
$\left\{ i_{1},\ldots,i_{m_{1}+n_{2}+m_{2}+n_{1}-1}\right\} $ of
interfaces by two empty interfaces, $i_{0}$ and $i_{m_{1}+n_{2}+m_{2}+n_{1}}$,
located, respectively, above and below $\boldsymbol{B}$.

Let us further refer to all the edges and vertices of $\boldsymbol{B}$
located (at least partially) between the elementary cuts $I_{k-1}$
and $I_{k}$ as the \emph{fragment} $f_{k}$ of $\boldsymbol{B}$.
We augment the set of fragments with two additional fragments: $f_{1}$
including all edges and vertices located (at least partially) above
the elementary cut $I_{1}$ and $f_{m_{1}+n_{2}+m_{2}+n_{1}}$ including
all edges and vertices located (at least partially) below the elementary
cut $I_{m_{1}+n_{2}+m_{2}+n_{1}-1}$. It is clear from these definitions
that for $1\leq k\leq m_{1}+n_{2}+m_{2}+n_{1}$
\begin{itemize}
\item $i_{k-1}\subset f_{k}$ is the upper interface of $f_{k}$,
\item $i_{k}\subset f_{k}$ is the lower interface of $f_{k}$. 
\end{itemize}
Consider now a fragment $f_{k}\subset\boldsymbol{B}$ together with
its upper interface $i_{k-1}$ and its lower interface $i_{k}$. Denote
by $e_{\text{first}}$ the left-most vertical edge of $f_{k}$, and
by $e_{\text{last}}$, the right-most vertical edge of $f_{k}$. We
can now define the function $\text{shape}$ as follows
\begin{equation}
\text{shape}\left(f_{k}\right)=\begin{cases}
\mathtt{W} & \text{when }e_{\text{first}}\in i_{k}\text{ and }e_{\text{last}}\in i_{k}\\
\mathtt{N} & \text{when \ensuremath{e_{\text{first}}\in i_{k-1}} and \ensuremath{e_{\text{last}}\in i_{k-1}}}\\
\mathtt{R} & \text{when \ensuremath{e_{\text{first}}\in i_{k-1}} and \ensuremath{e_{\text{last}}\in i_{k}}}\\
\mathtt{L} & \text{when }e_{\text{first}}\in i_{k}\text{ and }e_{\text{last}}\in i_{k-1}
\end{cases}
\end{equation}
The symbols $\mathtt{W}$, $\mathtt{N}$, $\mathtt{R}$, and $\mathtt{L}$
describe geometrically the shape (respectively: wider, narrower, to-the-right,
and to-the-left) of each fragment. 
 Using this terminology, it is possible to apply the function $\text{shape}$
to $\boldsymbol{B}=\left(f_{1},f_{2},\ldots,f_{m_{1}+n_{2}+m_{2}+n_{1}}\right)$,
simply by mapping it to the sequence of its fragments. For example,
the shape of the structure shown in Fig.~\ref{fig:graphdef} is specified
by the following sequence
\begin{equation}
\text{shape}\left(Rb\left(3,6,5,4\right)\right)=\mathtt{WWWWWLLLLNNRRRRNNN}
\end{equation}

Let us now consider an arbitrary Clar cover $\boldsymbol{C}$ of $\boldsymbol{B}$.
For every edge $e$ of $\boldsymbol{B}$ , we define a \emph{covering
order} function $\text{ord}\left(e\right)$ as follows
\begin{equation}
\text{ord}\left(e\right)=\begin{cases}
1 & \text{when \ensuremath{\exists\,K_{2}\subset\boldsymbol{C}:e\in K_{2}}}\\
\frac{1}{2} & \text{when }\ensuremath{\exists\,C_{6}\subset\boldsymbol{C}:e\in C_{6}}\\
0 & \text{otherwise}
\end{cases}
\end{equation}

\noindent This definition can be naturally extended to \emph{covering
order} \emph{of interfaces} by defining the order of the interface
$\text{ord}\left(i\right)$ as
\begin{equation}
\text{ord}\left(i\right)=\sum_{e\in i}\text{ord}\left(e\right).
\end{equation}
The interfaces $i_{0}$ and $i_{m_{1}+n_{2}+m_{2}+n_{1}}$are empty,
thus
\begin{equation}
\text{ord}\left(i_{0}\right)=0=\text{ord}\left(i_{m_{1}+n_{2}+m_{2}+n_{1}}\right).
\end{equation}
 It turns out that the orders of the remaining interfaces can be conveniently
computed in an iterative fashion using 
\begin{thm}
\textbf{\emph{\label{thm:1st rule-1}(First rule of interface theory:
interface order criterion) }}\citep{langner2019IFTTheorems5,langner2019BasicApplications6}
Let $\boldsymbol{C}$ be some Clar cover of a benzenoid $\boldsymbol{B}$.
Let $f_{k}$ be a fragment of $\boldsymbol{B}$, and let $i_{k-1}$
and $i_{k}$ be the upper and lower interfaces of $f_{k}$, respectively.
The following conditions are always satisfied.
\begin{lyxlist}{00.00.0000}
\item [{$(a)$}] If $f_{k}$ has the shape $\mathtt{W}$, then \emph{$\text{ord}(i_{k})=\text{ord}(i_{k-1})+1$.}
\item [{$(b)$}] If $f_{k}$ has the shape $\mathtt{N}$, then \emph{$\text{ord}(i_{k})=\text{ord}(i_{k-1})-1$.}
\item [{$(c)$}] If $f_{k}$ has the shape $\mathtt{R}$ or $\mathtt{L}$
, then \emph{$\text{ord}(i_{k})=\text{ord}(i_{k-1})$.}
\end{lyxlist}
\end{thm}

\noindent The interface orders determined this way starting with $\text{ord}\left(i_{0}\right)=0$
depend only on the shape of $\boldsymbol{B}$ and are independent
of the choice of $\boldsymbol{C}$. Therefore, the interface orders
are identical for every Clar cover $\boldsymbol{C}$. It is straightforward
to show that the interface orders of $Rb\left(3,6,5,4\right)$ shown
in Fig.~\ref{fig:graphdef} are specified by the following sequence
\[
\left(\text{ord}\left(i_{0}\right),\ldots,\text{ord}\left(i_{m_{1}+n_{2}+m_{2}+n_{1}}\right)\right)=\left(0,1,2,3,4,5,5,5,5,5,4,3,3,3,3,3,2,1,0\right)
\]
and the interface orders of $Rb\left(5,9,4,6\right)$ shown in Fig.~\ref{fig:interfaces}
are specified by the following sequence
\[
\left(\text{ord}\left(i_{0}\right),\ldots,\text{ord}\left(i_{m_{1}+n_{2}+m_{2}+n_{1}}\right)\right)=\left(0,1,2,3,4,\underbrace{5}_{10\text{ times}},\underbrace{6}_{4\text{ times}},5,4,3,2,1,0\right)
\]

Let us now consider in detail the interface $i_{m_{1}+n_{2}}$ of
$Rb\left(n_{1},n_{2},m_{1},m_{2}\right)$. Application of Theorem~\ref{thm:1st rule-1}
to $i_{m_{1}+n_{2}}$ shows that $\text{ord}(i_{m_{1}+n_{2}})=N\equiv\min\left(m_{1},n_{1}\right)$.
At the same time simple geometrical considerations show that $i_{m_{1}+n_{2}}$
consists of $N+1$ vertical edges $e_{0},\ldots,e_{N}$, where the
numbering proceeds from right to left. Since 
\begin{equation}
\text{ord}(i_{m_{1}+n_{2}})=\sum_{k=0}^{N}\text{ord}\left(e_{k}\right)=N\label{eq:m1n2ord}
\end{equation}
and since each $\text{ord}\left(e_{k}\right)$ can take on only three
values: $0$, $\frac{1}{2}$, and $1$, the interface order $\text{ord}(i_{m_{1}+n_{2}})=N$
can be created from the interface edges orders $\text{ord}\left(e_{k}\right)$
only in two possible ways
\begin{eqnarray}
N & = & \underbrace{1+\ldots+1}_{N\text{ times}}~~~+~~~\underbrace{0}_{1\text{ time}}\label{eq:sing}\\
N & = & \underbrace{1+\ldots+1}_{N-1\text{ times}}~~~+~~~\underbrace{{\textstyle \frac{1}{2}}+{\textstyle \frac{1}{2}}}_{2\text{ times}}\label{eq:ring}
\end{eqnarray}
The first of these two choices, described by Eq.~(\ref{eq:sing}),
correspond to Clar covers in which the interface $i_{m_{1}+n_{2}}$
is composed of $N$ double bonds and $1$ single bond. Clearly, there
exist $N+1$ distinct classes of Clar covers fulfilling this condition;
in each distinct class, the single bond is located at the position
$e_{k}$ with $k\in\left\{ 0,\ldots,N\right\} $. Each of such distinct
classes corresponds to a single summand in Eq.~(\ref{eq:sum1}).
The second of these choices, described by Eq.~(\ref{eq:ring}), correspond
to Clar covers in which the interface $i_{m_{1}+n_{2}}$ is composed
of $N-1$ double bonds and $2$ aromatic bonds belonging to the same
hexagon $C_{6}$. Clearly, there exist $N$ distinct classes of Clar
covers fulfilling this condition; in each distinct class, the hexagon
$C_{6}$ is located at the positions $e_{k-1}$ and $e_{k}$ with
$k\in\left\{ 1,\ldots,N\right\} $. Each of such distinct classes
corresponds to a single summand in Eq.~(\ref{eq:sum2}). To complete
the proof of Eq.~(\ref{eq:ZZRbM}) it remains to be demonstrated
that for each of these distinct classes of Clar covers, the covering
orders of edges in the interface $i_{m_{1}+n_{2}}$ induce a fixed-bond
region in $\boldsymbol{B}$ separating two not-fixed-bond regions,
each of them in a shape of a parallelogram.

\begin{figure}
\centering{}\includegraphics[clip,scale=0.4]{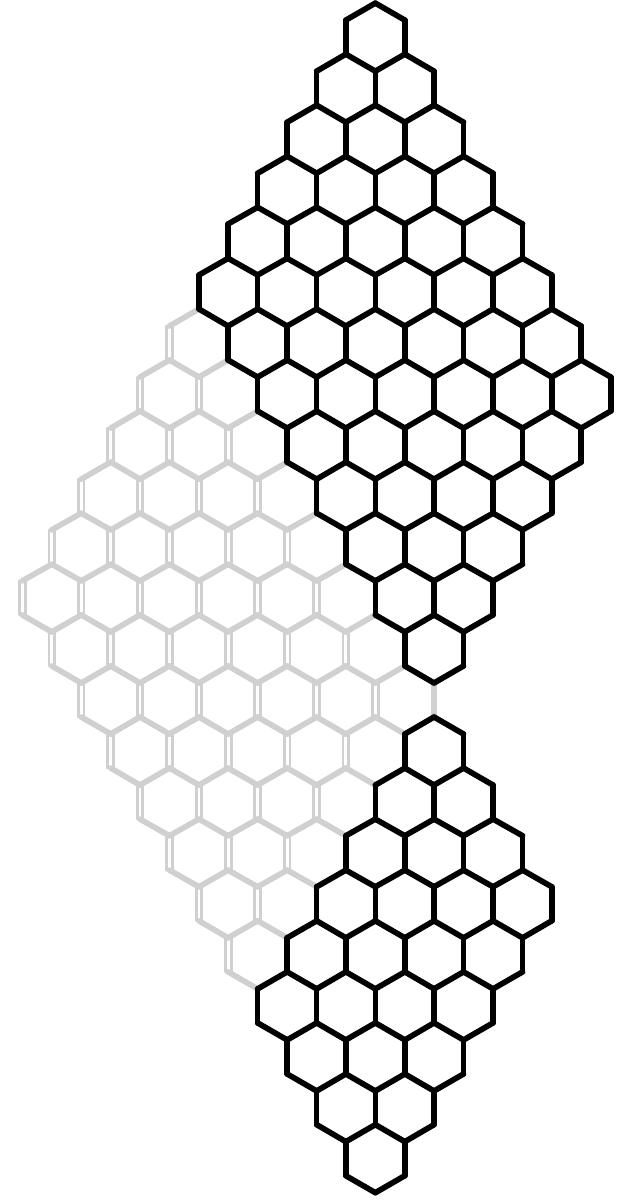}\hspace{2cm}\includegraphics[clip,scale=0.4]{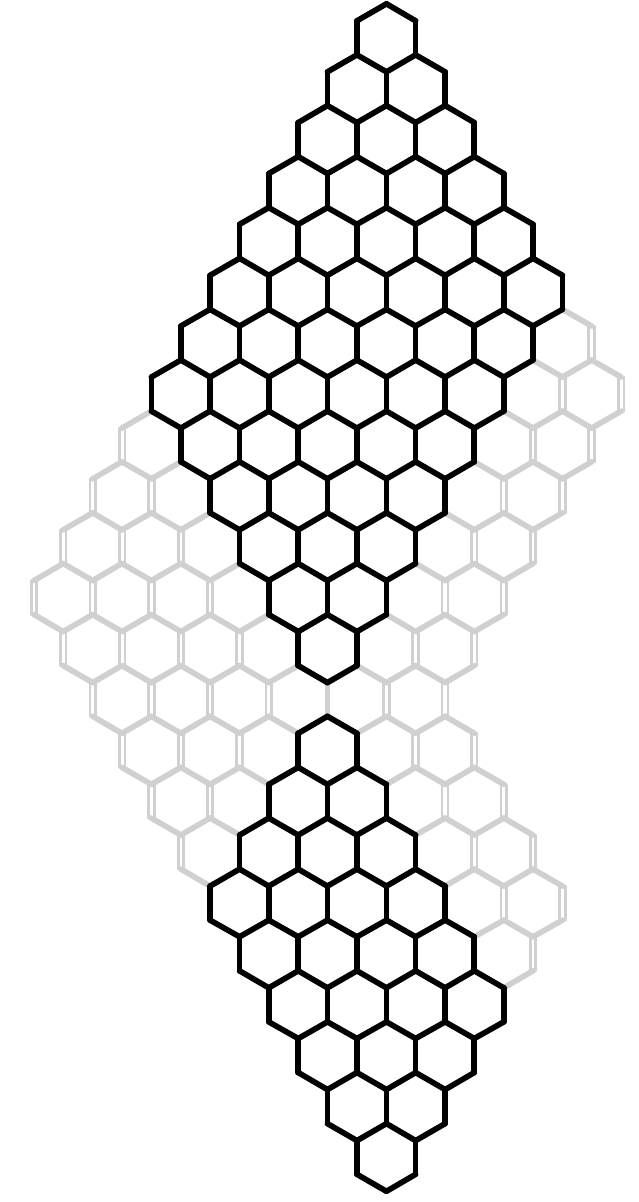}\hspace{2cm}\includegraphics[clip,scale=0.4]{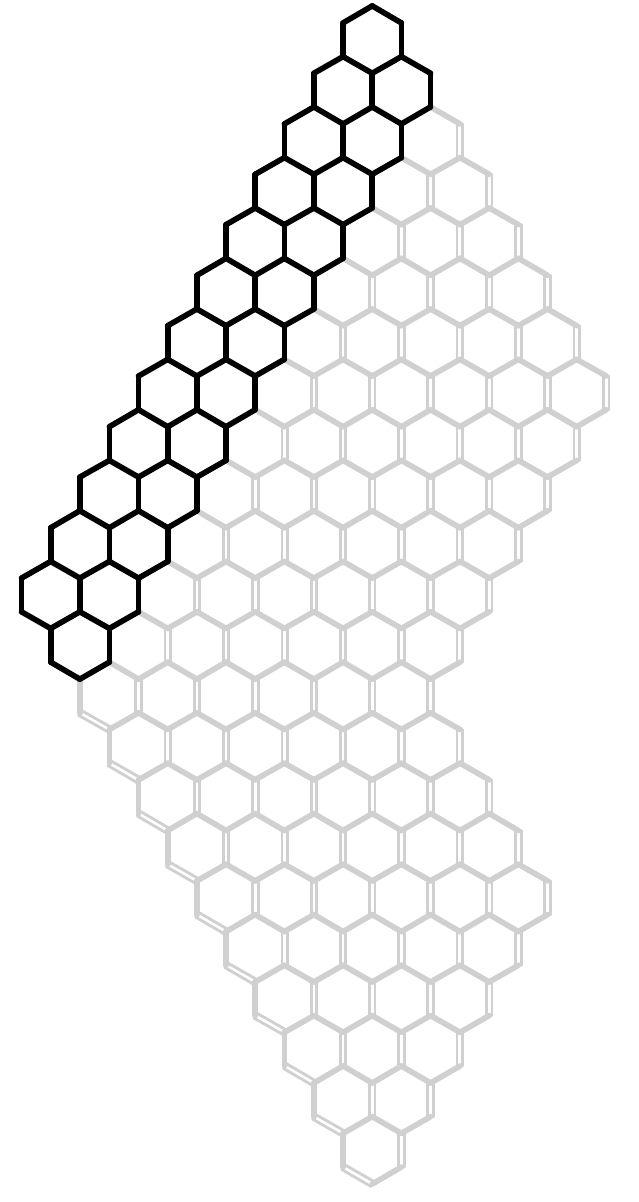}\caption{The covering of the interface $i_{m_{1}+n_{2}}$ involving $N$ double
bonds and $1$ single bond (in position $e_{0}$, $e_{k}$, and $e_{N}$,
respectively) induces a fixed-bond region (in gray) in $Rb\left(n_{1},n_{2},m_{1},m_{2}\right)$,
which separates two not-fixed-bond regions (in black), each of them
in a shape of a parallelogram. The position of the single bond determines
the shapes of the parallelograms. Here, $m_{1}=8$, $m_{2}=4$, $n_{1}=6$,
and $n_{2}=6$.\label{fig:rb1dec-1}}
\end{figure}

Let us first consider a class of Clar covers of $\boldsymbol{B}$
corresponding to Eq.~(\ref{eq:sing}) with a single bond in the position
$e_{k}$ of the interface $i_{m_{1}+n_{2}}$ and double bonds in its
remaining positions. Fig.~\ref{fig:rb1dec-1} 
 shows that the systems of double bonds in the interface $i_{m_{1}+n_{2}}$
propagates down and up in $\boldsymbol{B}$, uniquely deciding the
covering orders for a large portion of this structure. Each Clar cover
belonging to this class shares this region of fixed bonds. However,
there remain two disconnected regions in $\boldsymbol{B}$, each in
the shape of a parallelogram, for which the covering characters are
not determined by the covering of the interface $i_{m_{1}+n_{2}}$.
The sizes of these two parallelograms are determined by the structural
parameters $n_{1}$, $n_{2}$, $m_{1}$, and $m_{2}$ and the location
$k$ of the single bond. It is easy to see that the upper parallelogram
is $M\left(m_{1}-k,n_{2}+k\right)$ and the lower one is $M\left(m_{2}+k,n_{1}-k\right)$.
The product of the ZZ polynomials of these two parallelograms
\[
\text{ZZ}\left(M\left(m_{1}-k,n_{2}+k\right),x\right)\cdot\text{ZZ}\left(M\left(m_{2}+k,n_{1}-k\right),x\right)
\]
describes the contribution of this class of Clar covers to the ZZ
polynomial of $\boldsymbol{B}$. Sum of these contributions for $k\in\left\{ 0,\ldots,N\right\} $
reproduces Eq.~(\ref{eq:sum1}).

\begin{figure}
\centering{}\includegraphics[clip,scale=0.4]{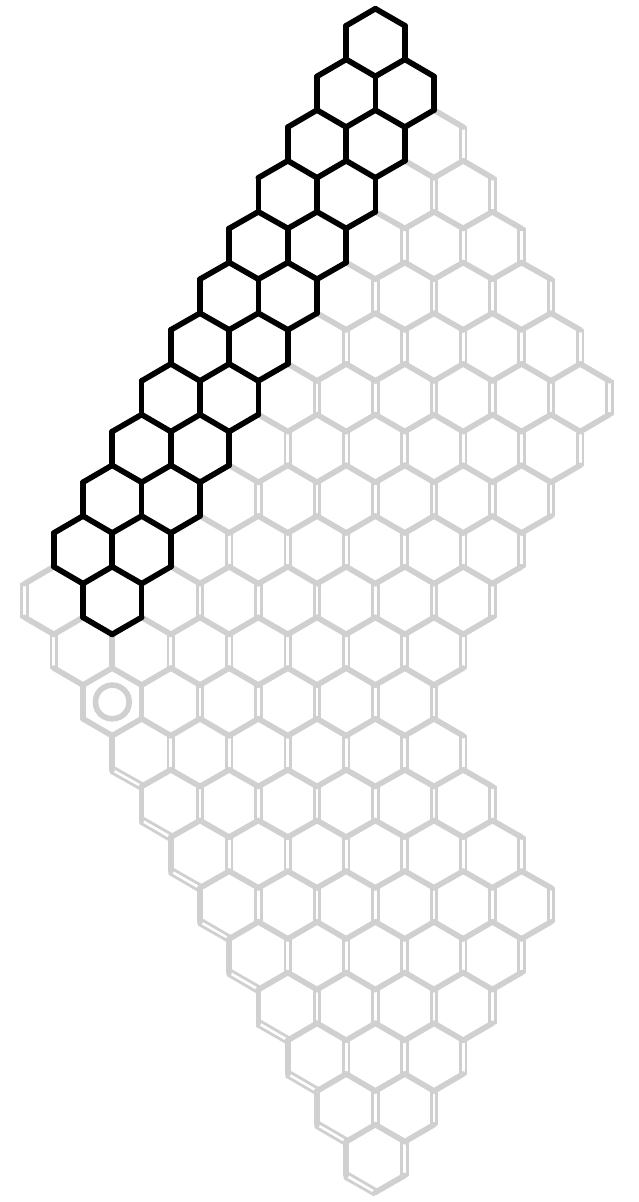}\hspace{2cm}\includegraphics[clip,scale=0.4]{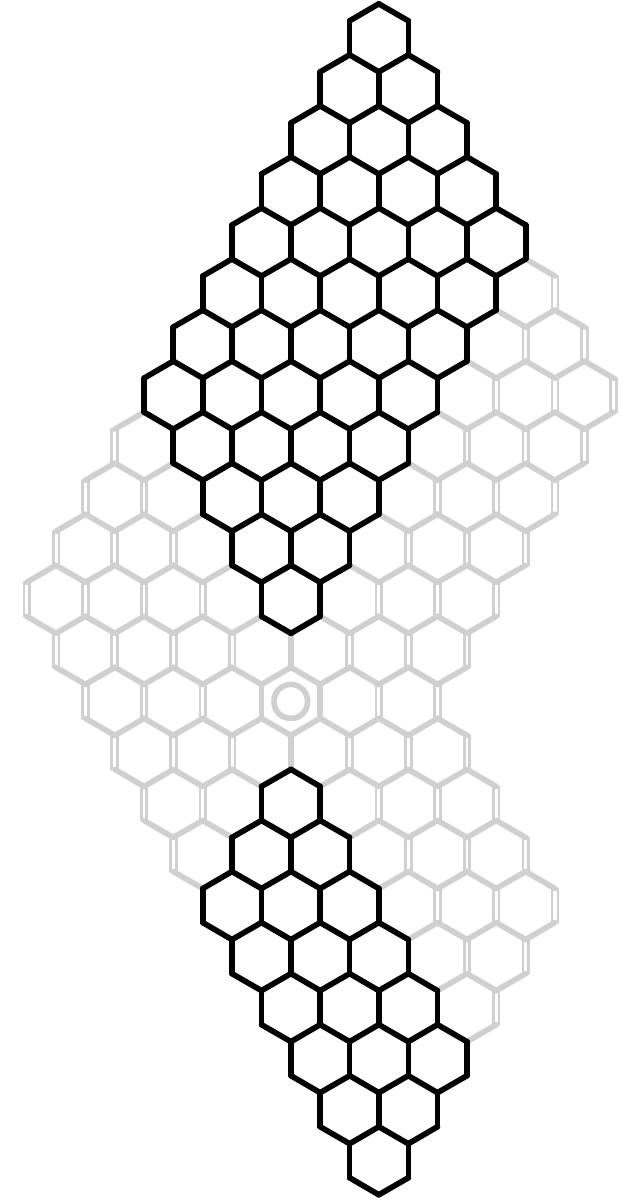}\hspace{2cm}\includegraphics[clip,scale=0.4]{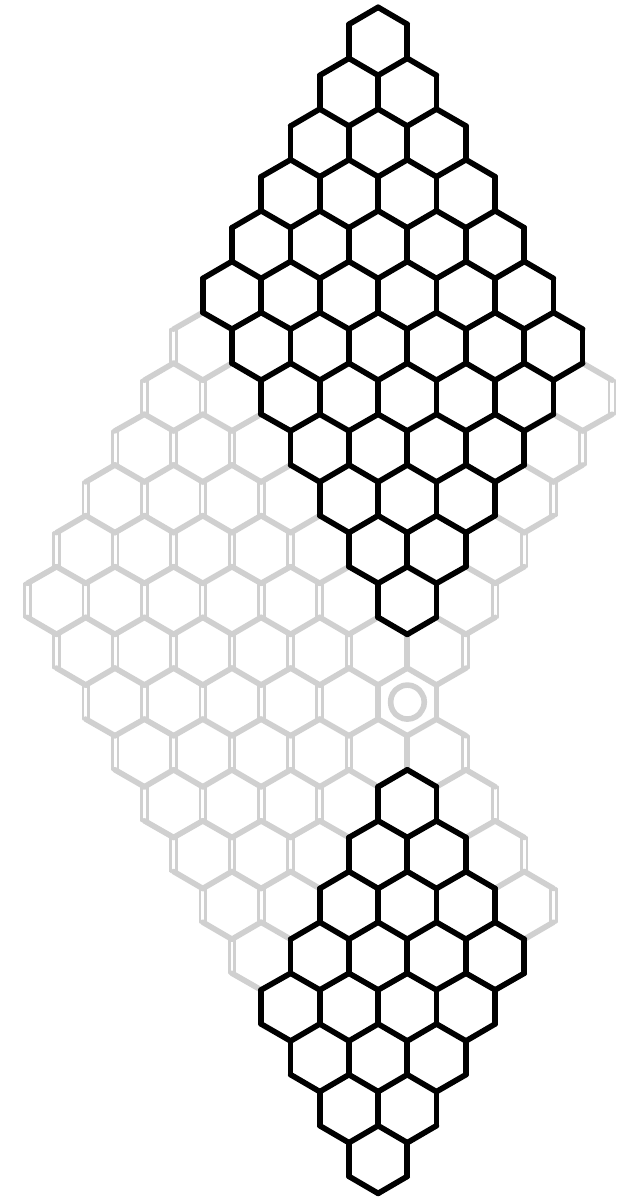}\caption{The covering of the interface $i_{m_{1}+n_{2}}$ involving $N-1$
double bonds and $1$ Clar sextet (in position $e_{1}$, $e_{k}$,
and $e_{N}$, respectively) induces a fixed-bond region (in gray)
in $Rb\left(n_{1},n_{2},m_{1},m_{2}\right)$, which separates two
not-fixed-bond regions (in black), each of them in a shape of a parallelogram.
The position of the Clar sextet determines the shapes of the parallelograms.
Here, $m_{1}=8$, $m_{2}=4$, $n_{1}=6$, and $n_{2}=6$.\label{fig:rb1dec-1-1}}
\end{figure}

Let us now consider a class of Clar covers of $\boldsymbol{B}$ corresponding
to Eq.~(\ref{eq:ring}) with a hexagon $C_{6}$ located at the positions
$e_{k-1}$ and $e_{k}$ of the interface $i_{m_{1}+n_{2}}$ and double
bonds in its remaining positions. Fig.~\ref{fig:rb1dec-1-1} 
 shows that the covered bonds in the interface $i_{m_{1}+n_{2}}$
induce a system of double bonds in the interfaces $i_{m_{1}+n_{2}-1}$
and $i_{m_{1}+n_{2}+1}$, which propagate down and up in $\boldsymbol{B}$,
again uniquely deciding the covering orders for a large portion of
this structure. Each Clar cover belonging to this class shares this
region of fixed bonds. Again, there remain two disconnected regions
in $\boldsymbol{B}$, each in the shape of a parallelogram, for which
the covering characters are not determined by the covering of the
interface $i_{m_{1}+n_{2}}$. The sizes of these two parallelograms
are determined by the structural parameters $n_{1}$, $n_{2}$, $m_{1}$,
and $m_{2}$ and the location $k$ of the Clar sextet. It is again
easy to see that the upper parallelogram is $M\left(m_{1}-k,n_{2}-1+k\right)$
and the lower one is $M\left(m_{2}-1+k,n_{1}-k\right)$. The product
of the ZZ polynomials of these two parallelograms
\[
\text{ZZ}\left(M\left(m_{1}-k,n_{2}-1+k\right),x\right)\cdot\text{ZZ}\left(M\left(m_{2}-1+k,n_{1}-k\right),x\right)
\]
describes the contribution of this class of Clar covers to the ZZ
polynomial of $\boldsymbol{B}$. Sum of these contributions for $k\in\left\{ 1,\ldots,N\right\} $
reproduces Eq.~(\ref{eq:sum2}) and concludes the proof of Eq.~(\ref{eq:ZZRbM}).

\section{Discussion and conclusion}

We have derived a closed-form formula for the ZZ polynomial of ribbons
$\boldsymbol{B}\equiv Rb\left(n_{1},n_{2},m_{1},m_{2}\right)$, an
important class of elementary pericondensed benzenoids. The formal
demonstration of its correctness is based on the recently developed
interface theory of benzenoids. The discovered formula
\[
\text{ZZ}\left(\boldsymbol{B},x\right)=~~~~~~~~~~~~~~~~~~~~~~~~~~~~~~~~~~~~~~~~~~~~~~~~~~~~~~~~~~~~~~~~~~~~
\]
\[
=~\sum_{k=0}^{\min\left(n_{1},m_{1}\right)}\sum_{i=0}^{m_{1}-k}\sum_{j=0}^{n_{1}-k}\binom{m_{1}-k}{j}\binom{n_{2}+k}{j}\binom{m_{2}+k}{i}\binom{n_{1}-k}{i}\left(1+x\right)^{i+j}+
\]
\[
~~~~~~~~+~x\sum_{k=1}^{\min\left(n_{1},m_{1}\right)}\sum_{i=0}^{m_{1}-k}\sum_{j=0}^{n_{1}-k}\binom{m_{1}-k}{j}\binom{n_{2}-1+k}{j}\binom{m_{2}-1+k}{i}\binom{n_{1}-k}{i}\left(1+x\right)^{i+j}
\]
uniquely determines the most important topological invariants of $Rb\left(n_{1},n_{2},m_{1},m_{2}\right)$: 
\begin{itemize}
\item the number of Kekulé structures
\[
K\left\{ \boldsymbol{B}\right\} =\sum_{k=0}^{\min\left(n_{1},m_{1}\right)}\sum_{i=0}^{m_{1}-k}\sum_{j=0}^{n_{1}-k}\binom{m_{1}-k}{j}\binom{n_{2}+k}{j}\binom{m_{2}+k}{i}\binom{n_{1}-k}{i}
\]
\item the number of Clar covers
\[
C\left\{ \boldsymbol{B}\right\} =\sum_{k=0}^{\min\left(n_{1},m_{1}\right)}\sum_{i=0}^{m_{1}-k}\sum_{j=0}^{n_{1}-k}\binom{m_{1}-k}{j}\binom{n_{2}+k}{j}\binom{m_{2}+k}{i}\binom{n_{1}-k}{i}2^{i+j}
\]
\[
+\sum_{k=1}^{\min\left(n_{1},m_{1}\right)}\sum_{i=0}^{m_{1}-k}\sum_{j=0}^{n_{1}-k}\binom{m_{1}-k}{j}\binom{n_{2}-1+k}{j}\binom{m_{2}-1+k}{i}\binom{n_{1}-k}{i}2^{i+j}
\]
\item the Clar number $Cl=\deg\left(\text{ZZ}\left(\boldsymbol{B},x\right)\right)$, 
\item and the number of Clar structures equal to $\text{coeff }\left(\text{ZZ}\left(\boldsymbol{B},x\right),x^{Cl}\right)$.
\end{itemize}
\noindent Interestingly, it is straightforward to obtain the Clar
number of $Rb\left(n_{1},n_{2},m_{1},m_{2}\right)$ and the number
of Clar structures of $Rb\left(n_{1},n_{2},m_{1},m_{2}\right)$ directly
from the ZZ polynomial, but it seems to be a formidable task to extract
these two quantities directly from the structural constants $n_{1}$,
$n_{2}$, $m_{1}$, and $m_{2}$. For example, at the moment, the
most compact formula for $Cl$ in terms of the structural constants
$n_{1}$, $n_{2}$, $m_{1}$, and $m_{2}$ that we are aware of is
given by the following expression
\begin{equation}
Cl=\max\left(Cl_{\text{s}},Cl_{\text{r}}\right),\label{eq:Cl}
\end{equation}
where
\begin{eqnarray}
Cl_{\text{s}} & = & \max_{k\in\left\{ 0,\ldots,\min\left(m_{1},n_{1}\right)\right\} }\left(\min\left(m_{1}-k,n_{2}+k\right)+\min\left(n_{1}-k,m_{2}+k\right)\right)\label{eq:Cls}\\
Cl_{\text{r}} & = & \max_{k\in\left\{ 1,\ldots,\min\left(m_{1},n_{1}\right)\right\} }\left(1+\min\left(m_{1}-k,n_{2}-1+k\right)+\min\left(n_{1}-k,m_{2}-1+k\right)\right)\label{eq:Clr}
\end{eqnarray}
Both of these terms are needed, as the following examples show. For
$Rb\left(2,2,1,1\right)$, we have $Cl_{\text{s}}=3$ and $Cl_{\text{r}}=2$,
so $Cl=Cl_{\text{s}}.$ For $Rb\left(3,2,1,2\right)$, we have $Cl_{\text{s}}=3$
and $Cl_{\text{r}}=4$, so $Cl=Cl_{\text{r}}.$ We believe that Eqs.~(\ref{eq:Cl})\textendash (\ref{eq:Clr})
cannot be simplified much. Therefore, it should be very instructive
to see this expression for researchers who try to determine Clar numbers
of simple benzenoids directly from geometrical considerations. The
structural complexity of this formula suggests that transforming relatively
easy geometrical constructs into an algebraic expression can be cumbersome.

The last two classes of elementary pericondensed benzenoids, for which
closed-form ZZ polynomial formulas remain to be found, are hexagonal
flakes $O\left(k,m,n\right)$ and oblate rectangles $Or\left(m,n\right)$.
We hope that our results will stimulate mathematicians and mathematically-oriented
chemists to discover these two last missing formulas.

\end{document}